\documentclass[12pt,twoside,leqno]{article}
\usepackage{amsmath}
\usepackage{amssymb}
\usepackage{amsxtra}
\usepackage{amscd}
\usepackage{amsthm}
\usepackage[mathscr]{eucal}
\usepackage{color}
\usepackage{fancyhdr}

\theoremstyle{plain}
\newtheorem{thm}[subsection]{Theorem}
\newtheorem{prop}[subsection]{Proposition}

\newtheorem{lem}[subsection]{Lemma}
\newtheorem{conj}[subsection]{Conjecture}

\theoremstyle{definition}

\newtheorem{rem}[subsection]{Remark}
\newtheorem{para}[subsection]{}

\newenvironment{pf}{\proof[\proofname]}{\endproof}

\begin{document}

\title{Logarithmic Tate conjectures over finite fields}

\author{Kazuya Kato, Chikara Nakayama, Sampei Usui}
\pagestyle{fancy}
\lhead{}
\rhead{{\sc Logarithmic Tate conjectures over finite fields}}
\renewcommand{\headrulewidth}{0pt}

\maketitle

\newcommand\Cal{\mathcal}
\newcommand\define{\newcommand}

\define\gp{\mathrm{gp}}%
\define\fs{\mathrm{fs}}%
\define\an{\mathrm{an}}%
\define\mult{\mathrm{mult}}%
\define\Ker{\mathrm{Ker}\,}%
\define\Coker{\mathrm{Coker}\,}%
\define\Hom{\mathrm{Hom}\,}%
\define\Ext{\mathrm{Ext}\,}%
\define\rank{\mathrm{rank}\,}%
\define\gr{\mathrm{gr}}%
\define\cHom{\Cal{Hom}}
\define\cExt{\Cal Ext\,}%

\define\cB{\Cal B}
\define\cC{\Cal C}
\define\cD{\Cal D}
\define\cO{\Cal O}
\define\cS{\Cal S}
\define\cT{\Cal T}
\define\cM{\Cal M}
\define\cG{\Cal G}
\define\cH{\Cal H}
\define\cE{\Cal E}
\define\cF{\Cal F}
\define\cN{\Cal N}
\define\cQ{\Cal Q}
\define\fF{\frak F}
\define{\sW}{\Cal W}

\define\Dc{\check{D}}
\define\Ec{\check{E}}

\newcommand{\N}{{\mathbb{N}}}
\newcommand{\Q}{{\mathbb{Q}}}
\newcommand{\Z}{{\mathbb{Z}}}
\newcommand{\R}{{\mathbb{R}}}
\newcommand{\C}{{\mathbb{C}}}
\newcommand{\bN}{{\mathbb{N}}}
\newcommand{\bQ}{{\mathbb{Q}}}
\newcommand{\bF}{{\mathbb{F}}}
\newcommand{\bZ}{{\mathbb{Z}}}
\newcommand{\bA}{{\mathbb{A}}}
\newcommand{\bP}{{\mathbb{P}}}
\newcommand{\bR}{{\mathbb{R}}}
\newcommand{\bC}{{\mathbb{C}}}
\newcommand{\bG}{{\mathbb{G}}}
\newcommand{\bbQ}{{\bar \mathbb{Q}}}
\newcommand{\ol}[1]{\overline{#1}}
\newcommand{\too}{\longrightarrow}
\newcommand{\respect}{\rightsquigarrow}
\newcommand{\compatible}{\leftrightsquigarrow}
\newcommand{\upc}[1]{\overset {\lower 0.3ex \hbox{${\;}_{\circ}$}}{#1}}
\newcommand{\Gmlog}{\bG_{m, \log}}
\newcommand{\Gm}{\bG_m}
\newcommand{\ep}{\varepsilon}
\newcommand{\Spec}{\operatorname{Spec}}
\newcommand{\nilp}{\operatorname{nilp}}
\newcommand{\prim}{\operatorname{prim}}
\newcommand{\val}{{\mathrm{val}}} 
\newcommand{\n}{\operatorname{naive}}
\newcommand{\bs}{\operatorname{\backslash}}
\newcommand{\Gal}{\operatorname{{Gal}}}
\newcommand{\gal}{{\rm {Gal}}({\bar \Q}/{\Q})}
\newcommand{\galp}{{\rm {Gal}}({\bar \Q}_p/{\Q}_p)}
\newcommand{\gall}{{\rm{Gal}}({\bar \Q}_\ell/\Q_\ell)}
\newcommand{\wep}{W({\bar \Q}_p/\Q_p)}
\newcommand{\wel}{W({\bar \Q}_\ell/\Q_\ell)}
\newcommand{\Ad}{{\rm{Ad}}}
\newcommand{\BS}{{\rm {BS}}}
\newcommand{\even}{\operatorname{even}}
\newcommand{\End}{{\rm {End}}}
\newcommand{\odd}{\operatorname{odd}}
\newcommand{\GL}{\operatorname{GL}}
\newcommand{\np}{\text{non-$p$}}
\newcommand{\g}{{\gamma}}
\newcommand{\G}{{\Gamma}}
\newcommand{\Lam}{{\Lambda}}
\newcommand{\La}{{\Lambda}}
\newcommand{\lam}{{\lambda}}
\newcommand{\la}{{\lambda}}
\newcommand{\uL}{{{\hat {L}}^{\rm {ur}}}}
\newcommand{\uQp}{{{\hat \Q}_p}^{\text{ur}}}
\newcommand{\sel}{\operatorname{Sel}}
\newcommand{\dt}{{\rm{Det}}}
\newcommand{\Sig}{\Sigma}
\newcommand{\fil}{{\rm{fil}}}
\newcommand{\SL}{{\rm{SL}}}
\renewcommand{\sl}{{\frak{sl}}}%
\newcommand{\spl}{{\rm{spl}}}
\newcommand{\st}{{\rm{st}}}
\newcommand{\Isom}{{\rm {Isom}}}
\newcommand{\Mor}{{\rm {Mor}}}
\newcommand{\bg}{\bar{g}}
\newcommand{\id}{{\rm {id}}}
\newcommand{\cone}{{\rm {cone}}}
\newcommand{\al}{a}
\newcommand{\ChL}{{\cal{C}}(\La)}
\newcommand{\Image}{{\rm {Image}}}
\newcommand{\toric}{{\operatorname{toric}}}
\newcommand{\torus}{{\operatorname{torus}}}
\newcommand{\Aut}{{\rm {Aut}}}
\newcommand{\Qp}{{\mathbb{Q}}_p}
\newcommand{\barQp}{{\mathbb{Q}}_p}
\newcommand{\Qpur}{{\mathbb{Q}}_p^{\rm {ur}}}
\newcommand{\Zp}{{\mathbb{Z}}_p}
\newcommand{\Zl}{{\mathbb{Z}}_l}
\newcommand{\Ql}{{\mathbb{Q}}_l}
\newcommand{\Qlur}{{\mathbb{Q}}_l^{\rm {ur}}}
\newcommand{\F}{{\mathbb{F}}}
\newcommand{\eps}{{\epsilon}}
\newcommand{\epsLa}{{\epsilon}_{\La}}
\newcommand{\epsLaVxi}{{\epsilon}_{\La}(V, \xi)}
\newcommand{\epsOLaVxi}{{\epsilon}_{0,\La}(V, \xi)}
\newcommand{\Qplin}{{\mathbb{Q}}_p(\mu_{l^{\infty}})}
\newcommand{\otimesQplin}{\otimes_{\Qp}{\mathbb{Q}}_p(\mu_{l^{\infty}})}
\newcommand{\galFl}{{\rm{Gal}}({\bar {\Bbb F}}_\ell/{\Bbb F}_\ell)}
\newcommand{\gallur}{{\rm{Gal}}({\bar \Q}_\ell/\Q_\ell^{\rm {ur}})}
\newcommand{\galFF}{{\rm {Gal}}(F_{\infty}/F)}
\newcommand{\galFv}{{\rm {Gal}}(\bar{F}_v/F_v)}
\newcommand{\galF}{{\rm {Gal}}(\bar{F}/F)}
\newcommand{\epsVxi}{{\epsilon}(V, \xi)}
\newcommand{\epsOVxi}{{\epsilon}_0(V, \xi)}
\newcommand{\plim}{\lim_
{\scriptstyle 
\longleftarrow \atop \scriptstyle n}}
\newcommand{\sig}{{\sigma}}
\newcommand{\ga}{{\gamma}}
\newcommand{\del}{{\delta}}
\newcommand{\Vss}{V^{\rm {ss}}}
\newcommand{\Bst}{B_{\rm {st}}}
\newcommand{\Dpst}{D_{\rm {pst}}}
\newcommand{\Dcrys}{D_{\rm {crys}}}
\newcommand{\DdR}{D_{\rm {dR}}}
\newcommand{\Fin}{F_{\infty}}
\newcommand{\Kla}{K_{\lambda}}
\newcommand{\Ola}{O_{\lambda}}
\newcommand{\Mla}{M_{\lambda}}
\newcommand{\Det}{{\rm{Det}}}
\newcommand{\Sym}{{\rm{Sym}}}
\newcommand{\LaSa}{{\La_{S^*}}}
\newcommand{\cX}{{\cal {X}}}
\newcommand{\MHG}{{\frak {M}}_H(G)}
\newcommand{\tauMla}{\tau(M_{\lambda})}
\newcommand{\Fvur}{{F_v^{\rm {ur}}}}
\newcommand{\Lie}{{\rm {Lie}}}
\newcommand{\cL}{{\cal {L}}}
\newcommand{\cW}{{\cal {W}}}
\newcommand{\fq}{{\frak {q}}}
\newcommand{\cont}{{\rm {cont}}}
\newcommand{\SC}{{SC}}
\newcommand{\Om}{{\Omega}}
\newcommand{\dR}{{\rm {dR}}}
\newcommand{\crys}{{\rm {crys}}}
\newcommand{\hatSig}{{\hat{\Sigma}}}
\newcommand{\rdet}{{{\rm {det}}}}
\newcommand{\ord}{{{\rm {ord}}}}
\newcommand{\BdR}{{B_{\rm {dR}}}}
\newcommand{\BdRO}{{B^0_{\rm {dR}}}}
\newcommand{\Bcrys}{{B_{\rm {crys}}}}
\newcommand{\Qw}{{\mathbb{Q}}_w}
\newcommand{\barkappa}{{\bar{\kappa}}}
\newcommand{\cP}{{\Cal {P}}}
\newcommand{\cZ}{{\Cal {Z}}}
\newcommand{\cA}{{\Cal {A}}}
\newcommand{\oppLa}{{\Lambda^{\circ}}}
\newcommand{\sat}{{{\rm sat}}}
\renewcommand{\bar}{\overline}
\newcommand{\et}{\mathrm{\acute{e}t}}
\newcommand{\loget}{\mathrm{log\acute{e}t}}
\newcommand{\pri}{{{\rm prim}}}
\newcommand{\add}{{\rm{add}}}
\newcommand{\zar}{{{\rm zar}}}
\define\cR{\Cal R}
\newcommand{\br}{{\bold r}}
\newcommand{\mild}{{{\rm {mild}}}}

\begin{abstract} We formulate an analogue of Tate conjecture on algebraic cycles, for the log geometry over a finite field. We show that the weight-monodromy conjecture follows from this conjecture and from the semi-simplicity of the Frobenius action.  This conjecture suggests the existence of the monodromy cycle which gives the monodromy operator and an action of $\sl(2)$ on the cohomology, and which lives in the world of log motives. 
\end{abstract}

\renewcommand{\thefootnote}{\fnsymbol{footnote}}



\section*{Contents}

\noindent Introduction 

\noindent \S\ref{s:Tate}. Log Tate conjecture

\noindent \S\ref{s:WM}. Weight-monodromy conjecture

\noindent \S\ref{s:logK}. Log $K_0$ and the finer log Tate conjecture

\noindent \S\ref{s:mod}. The monodromy cycle

\noindent \S\ref{s:Ex}. Examples

\noindent \S\ref{s:theta}. Theta functions

\noindent References

\section*{Introduction}\label{s:intro}

Let $V$ be a projective smooth algebraic variety over a $p$-adic local field $K$.  The monodromy operator acts on the $\ell$-adic \'etale cohomology group $H^m_{\et}(V_{\bar K}, \Q_{\ell})$ for each prime number $\ell\neq p$.

There are problems on the monodromy operator:

1. Weight-monodromy conjecture. In the case where the residue field is finite, this conjecture gives a strong relation of the weight of the action of the Frobenius and the action of the monodromy operator on the $\ell$-adic \'etale cohomology group.

2. The monodromy operator  often has properties which are independent of $\ell$. This independence suggests that the monodromy operator is in fact an algebraic cycle in the algebraic geometry over the residue field $k$ of $K$. This subject was studied in C.\ Consani (\cite{CB}) (with Appendix by S.\ Bloch), C.\ Consani and M. Kim \cite{CK}, and also in \cite{logmot} Appendix.

These problems can be regarded as problems in logarithmic  geometry over the residue field $k$. Assume that $V$ has a strictly semistable reduction $X$. Then $X$ is  regarded as a log scheme over the standard log point $s$ over $k$. We have the identification 
$$H^m_{\et}(V_{\bar K}, \Q_{\ell})=H^m_{\loget}(X_{{\bar s}(\log)}, \Q_{\ell})$$ with the log \'etale cohomology group (\cite{Nc}, \cite{Nc2}), where $\bar s(\log)$ is the log separable closure of $s$.
Furthermore,  the action of $\Gal(\bar K/K)$ on the  left-hand-side  is identified with the action of the log fundamental group $\pi_1^{\log}(s)=\Gal(K^{\text{tame}}/K)$ on the  right-hand-side, where $K^{\text{tame}}\subset \bar K$ denotes the maximal tame extension of $K$.

In this paper, assuming that $k$ is a finite field, we formulate logarithmic analogues of the Tate conjecture on algebraic cycles. Such a conjecture was formulated and discussed in \cite{IKNU} Section 6 and in \cite{logmot} Appendix without assuming the finiteness of $k$, but here in the case where $k$ is a finite field, we formulate a stronger  log Tate conjecture \ref{Tate4} using a log Tate curve over $s$. A log Tate curve is the reduction of a Tate elliptic curve over $K$, which is regarded as a log scheme over $s$. Our method  is based on the fact that all log Tate curves over $s$ are isogenous in the case where $k$ is finite.  Note that the above problems 1 and 2 may be reduced to the finite residue field case by specialization arguments. 

The log Tate conjecture \ref{Tate4} says that if $m\geq 2r$, the map
$$\Q_{\ell} \otimes \gr^{m-r}K_{0,\lim}(X \times E^{m-2r})\to H^m_{\loget}(X_{{\bar s}(\log)}, \Q_{\ell})(r)^G  $$ 
is surjective. Here $E$ is the log Tate curve over $s$, $X \times E^{m-2r}$ means the fiber product over $s$, $K_{0,\lim}$ is the inductive limit of the $K$-groups $K_0$ for the log modifications of $X\times E^{m-2r}$,  $\gr^{m-r}$  is for the $\gamma$-filtration, 
$G=\pi_1^{\log}(s)$, and $(\cdot)^G$ is the part fixed by $G$.

We prove that if this log Tate conjecture is true and a certain semi-simplicity of the Frobenius action is true, then the weight-monodromy conjecture is true. 

In Conjecture \ref{conjiso2}, we  formulate a finer version
$$\Q_{\ell} \otimes K(X, m, r) \overset{\cong}\to H^m_{\loget}(X_{{\bar s}(\log)}, \Q_{\ell})(r)^G $$
of the log Tate conjecture \ref{Tate4}, by using a subquotient  $K(X, m, r)$ of $\Q\otimes \gr^{m-r}K_{0,\lim}(X \times E^{m-2r})$. 

The  log Tate conjecture \ref{conjiso2} suggests that there is a unique element of $K(X\times X, 2d, d-1)$ with $d=\dim X$  which induces the monodromy operator on $H^m_{\et}(V_{\bar K}, \Q_{\ell})$ for every $m$ and  $\ell\neq p$ and gives an action of $\sl(2)$ on it.  
  We call this conjectural element the {\it monodromy cycle}. In \ref{try}, we give an interpretation of the monodromy cycle by using log motives.

The above all seem to tell that the log Tate curve over a finite field has a special importance in arithmetic geometry. 

  K.\ Kato was 
partially supported by NSF grants DMS 1601861 and DMS 2001182.
C.\ Nakayama was 
partially supported by JSPS Grants-in-Aid for Scientific Research (C) 16K05093 and (C) 21K03199.
S.\ Usui was 
partially supported by JSPS Grants-in-Aid for Scientific Research (C) 17K05200.

\section{Log Tate conjecture}\label{s:Tate}

\begin{para}
\label{setting1} Let $k$ be a finitely generated field over the prime field.

The usual Tate conjecture is that for a projective smooth algebraic variety $X$ over $k$, the \'etale cohomology $H^{2r}_{\et}(X\otimes_k \bar k, \Q_{\ell})(r)^G$,  where $G=\Gal(\overline k/k)$, is generated by the classes of algebraic cycles on $X$ of codimension $r$. It is known that $H^m_{\et}(X\otimes_k \bar k, \Q_{\ell})(r)^G=0$ unless $m=2r$. 

\end{para}

\begin{para}
\label{setting2} 

Consider a log scheme  $s=\Spec(k)$ which is endowed with a log structure $M_s$ such that $M_s/\cO_s^\times=\N$. 

By Hilbert 90, there is a section $\pi$ of $M_s$ such that $\N \overset{\cong}\to M_s/\cO_s^\times\;;\; 1\mapsto \pi$. Such a $\pi$ is called a {\it generator} of the log structure of $s$.

Let 
$\frak S=\frak S(s)$ be the category of projective vertical log smooth fs log schemes over $s$ which have charts of the log structure Zariski locally. The reason why we put the last condition on the log structure is explained in \ref{K0lim} below.

By the above remark on the existence of a section $\pi$, $s$ itself belongs to $\frak S$.

For example, if $K, k$, and $V$ with strictly semistable reduction $X$ are as in Introduction, that is, if there is a projective smooth regular flat scheme $\frak X$ over the valuation ring $O_K$ of $K$ of strictly semistable reduction and $X=\frak X \otimes_{O_K} k$, then we can take $\Spec(k)$ with the inverse image $M_s$ of the natural log structure of $\Spec(O_K)$ as our $s$, a prime element of $K$ gives a generator of the log structure of $s$, and  $X$ with the inverse image of the  natural  log structure of $\frak X$ belongs to $\frak S$.

 Throughout this paper, $X$ denotes an object of $\frak S$ unless otherwise stated.

For a prime number $\ell$ which is different from the characteristic of $k$, we denote 
$$H^m(X)_{\ell}:= H^m_{\loget}(X_{{\bar s}(\log)}, \Q_{\ell}).$$

\end{para}

\begin{para}

In this situation,  $H^m(X)_{\ell}(r)^G$, where $G=\pi_1^{\log}(s)$, need not be zero for $m\neq 2r$. 

The following conjecture is well-known at least if $X$ is the reduction of $V$ as in Introduction.
\end{para}

\begin{conj}\label{Tate1} 
 $H^m(X)_{\ell}(r)^G=0$ unless $m\geq 2r$. 
\end{conj}

Now we consider the case $m\geq 2r$. 

\begin{para}\label{K0lim}  Let
$$K_{0,\lim}(X)= \varinjlim_{X'}\; K_0(X'),$$
where $X'$ ranges over all log modifications (\cite{IKNU} 2.3.6) of $X$.

The condition that $X$ has charts of the log structure Zariski locally ensures that $X$ has sufficiently many 
log modifications to get a nice $K_{0,\lim}(X)$ as is explained in \cite{IKNU} Sections 2.2--2.4.

\end{para}

First, the case $m=2r$ is formulated in the similar style to the classical case.

\begin{conj}\label{Tate2} 

The Chern class map $$\Q_{\ell} \otimes \gr^rK_{0,\lim}(X)\to H^{2r}(X)_{\ell}(r)^G$$ is surjective.

\end{conj}

This  conjecture was discussed in \cite{IKNU} Section 6.

For general $m\geq 2r$, we have 

\begin{conj}\label{Tate3} 
  Let $m \ge 2r$. Then the following holds. 

$(1)$ The map 
$$\Q_{\ell} \otimes \gr^{m-r}K_{0,\lim}(X\times \mathbb{G}_m^{m-2r})\to H^m(X)_{\ell}(r)^G$$ is surjective.

$(2)$ The  Chern class map 
$$\Q_{\ell} \otimes \gr^rKH_{2r-m,\lim}(X)\to H^m(X)_{\ell}(r)^G$$ is surjective. Here $KH$ denotes the homotopy $K$-theory ({\rm \cite{W1}}) and $KH_{i,\lim}(X)$ is defined as  $\varinjlim_{X'} \;KH_i(X')$ in the same way as in $\ref{K0lim}$. 
\end{conj}

Note that the index $2r-m$ in the $K$-group in (2) is $\leq 0$. 

 (1) implies (2) as is explained in \ref{Tatec}.

 \begin{para}
 The conjectures (1) and (2) in Conjecture \ref{Tate3}  were considered in the resp.\ part of \cite{logmot} A.14. Actually, there in \cite{logmot}, we assumed that $X$ was saturated over $s$. We put this assumption just because we were following the analogy with the theory of log Hodge structures with unipotent local monodromy and because the action of the inertia subgroup $\pi_1^{\log}(\bar s)= \Gal(\bar s(\log)/\bar s)$ of $\pi_1^{\log}(s)= \Gal(\bar s(\log)/s)$ on $H^m(X)_{\ell}$ is unipotent for such an $X$ (see Proposition \ref{satunip}).
 
\end{para}

\begin{para}\label{Eqn}
Here we review the log Tate curve over $s$ from the viewpoint of  the theory of log abelian varieties (\cite{KKN}). This viewpoint makes some stories (for example, the story of the multiplication by $a\in \Z$ in  \ref{Kmr}) transparent. 

Let $q$ be a section of the log structure of our log point $s$  which does not belong to $k^\times$. Then we have a log elliptic curve $E^{(q)}=\mathbb{G}_{m,\log}^{(q)}/q^{\Z}$ which is an abelian group sheaf functor on the category of fs log schemes over $s$  for the (classical) \'etale topology. Here $\mathbb{G}_{m,\log}$ is the sheaf $T\mapsto \Gamma(T, M^{\gp}_T)$ and $(\cdot)^{(q)}$ means the part consisting of $t\in M^{\gp}_T$ such that locally, we have  $q^m|t|q^n$ for some $m, n\in \Z$ such that $m\leq n$, that is, 
$tq^{-m}$ and $q^nt^{-1}$ belong to $M_T$. 
  This $E^{(q)}$ is not  a representable functor, but it has big representable subfunctors $E^{(q,n)}=\mathbb{G}_{m,\log}^{(q,n)}/q^{\Z}$ for integers $n \geq 1$, called the models of $E^{(q)}$, where $\mathbb{G}_{m,\log}^{(q,n)}$ is the part of 
$\mathbb{G}_{m,\log}$ consisting of $t$ such that locally we have
 $q^r|t^n |q^{r+1}$ for some $r\in \Z$. 
  For a local field $K$ as in Introduction, for a lifting $\tilde q$ of $q$ to $K^\times$, $E^{(q,n)}$ is a reduction of the Tate elliptic curve over $K$ of $q$-invariant $\tilde q$. As a scheme, $E^{(q,n)}$ is an $n$-gon.

For $n\geq 2$, $E^{(q,n)}$ belongs to $\frak S$ ($E^{(q,1)}$ does not have charts Zariski locally due to  the self-intersection of the $1$-gon). 
In this paper, a log Tate curve $E$ means $E=E^{(q,n)}$ for some $q$ and for some  $n\geq 2$. 

$E^{(q,n)}$ is semistable if the image of $q$ in $\N$ under $M_s/\cO_s^\times \cong \N$ is $n$ and is saturated over $s$ if this image is a multiple of $n$.

 Assume that $k$ is a finite field. Then all the log elliptic curves $E^{(q)}$ are isogenous. 
   Hence,  $\Q\otimes K_{0,\lim}$ and the log \'etale cohomology involving $E=E^{(q,n)}$ are independent of the choices of $q$ and $n$.

\end{para}

\begin{conj}\label{Tate4} 
  Let $m \ge 2r$. 
Assume that $k$ is a finite field. Let $E$ be a log Tate curve over the standard log point over $k$. Then the map
$$\Q_{\ell} \otimes \gr^{m-r}K_{0,\lim}(X\times E^{m-2r})\to H^m(X)_{\ell}(r)^G$$ is surjective.

\end{conj}

Note that since all log Tate curves over $s$ are isogenous by the finiteness of $k$,  this conjecture 
does not depend on the choice of the log Tate curve $E$.

\begin{para}\label{Tatec} 

Assume $m\geq 2r$.

We have a commutative diagram
$$\begin{matrix} \gr^{m-r}K_{0,\lim}(X  \times E^{m-2r})& \to   & \gr^{m-r}K_{0,\lim}(X \times \mathbb{G}_m^{m-2r})&\to & \gr^rKH_{2r-m,\lim}(X) 
 \\
\downarrow &&&& \downarrow \\
\Hom_G(\Sym^{m-2r}H^1(E)_{\ell}, H^m(X)_{\ell}(r)) && \to &&
   H^m(X)_{\ell}(r)^G. \end{matrix}$$
The upper left arrow is the pullback by the inclusion $\mathbb{G}_m\subset E$, and the lower arrow  is the restriction to $e_1^{m-2r}$, where $e_1$ is the canonical base of $H^1(E)_{\ell}^G=\Q_{\ell}$. 

The upper right arrow is defined by using the Mayer--Vietoris sequence for 
$\bP^1= (\bA^1)\cup (\bA^1)^-$, $\bA^1\cap (\bA^1)^-=\mathbb{G}_m$, repeatedly.

The composite map $\gr^{m-r}K_{0,\lim}(X \times \mathbb{G}_m^{m-2r})\to H^m(X)_{\ell}(r)^G$ is defined also as
$$\gr^{m-r}K_{0,\lim}(X \times \mathbb{G}_m^{m-2r})\to H^{2m-2r}(X\times \mathbb{G}_m^{m-2r})_{\ell}(m-r)^G $$ $$\to (H^m(X)_{\ell} \otimes H^1(\mathbb{G}_m)_{\ell}^{\otimes (m-2r)})(m-r)^G= H^m(X)_{\ell}(r)^G$$
by using the K\"unneth decomposition. 

The right vertical arrow was used in \cite{logmot} in the studies of motives in log geometry.

By this diagram, Conjecture \ref{Tate4} implies Conjecture \ref{Tate3} (1) and Conjecture \ref{Tate3} (1) implies Conjecture \ref{Tate3} (2). 

If we admit Conjecture \ref{Tate2}, the left vertical arrow induces a surjection from $\Q_{\ell}\otimes$ the upper $K$-group, and hence Conjecture \ref{Tate4} becomes equivalent to the statement that the lower arrow in this diagram  is surjective.

\end{para}

\section{Weight-monodromy conjecture}\label{s:WM}

\begin{para} We review the monodromy operator and the weight-monodromy conjecture. 
 Let $q$ be a section of the log structure of $s$ which does not belong to $k^\times$. Let $I=\pi_1^{\log}(\bar s)$. We have a canonical homomorphism 
 $$a_q: I\to \Z_{\ell}(1)\;;\; \sig\mapsto (\sig(q^{1/\ell^n})/q^{1/\ell^n})_{n\geq 1}$$
 with finite cokernel.
 
 Let $X$ be an object of $\frak S$ and let $m\geq 0$. Then for some open subgroup $I'$ of $I$, the action of $I'$ on $H^m(X)_{\ell}$ is unipotent  and factors through $a_q: I'\to \Z_{\ell}(1)$. Taking an element $\sig$ of $I'$ such that $a_q(\sig)\neq 0$, we define the monodromy operator as
 $$N_q:=a_q(\sig)^{-1}\log(\sig)\; :\; H^m(X)_{\ell}\to H^m(X)_{\ell}(-1),$$
which is independent of such a $\sig$. 

 If $q'$ is another section of the log structure of $s$ which does not belong to $k^\times$ and if $c$ is the rational number such that the classes of $q'$ and $q^c$ in $\Q\otimes M^{\gp}_s/\cO_s^\times\cong \Q$ coincide, we have $N_{q'}=c^{-1}N_q$.  
 
 In the situation where the choice of $q$ is not important (as in the situation in the rest of this Section \ref{s:WM}), we denote $N_q$ simply as $N$. 
 \end{para}
 
 \begin{para}
 
  In the case where $k$ is a finite field, the weight-monodromy conjecture says the following. For the monodromy filtration  $W$  of the monodromy operator $N$ on $H^m(X)_{\ell}$, the 
  action of the geometric Frobenius on $W_r/W_{r-1}$ is of weight $m+r$ for every $r$.

  Strong results on the weight-monodromy conjecture were obtained by P.\ Scholze in \cite{Sc}, but it is still a conjecture.

\end{para}

\begin{thm}
\label{t:WM} Assume the log Tate conjecture $\ref{Tate4}$ in general (not only for $X$ and $m$ below). Assume the following {\rm (i)} and {\rm (ii)}. 

{\rm (i)} There is a non-degenerate pairing 

$$H^m(X)_{\ell} \times H^m(X)_{\ell} \to \Q_{\ell}(-m)$$
which is compatible with the action of $G$.

{\rm (ii)} The Frobenius action on $H^m(X)_{\ell}^{N=0}$ is semisimple.

 Then the weight-monodromy conjecture for $X$ is true.

\end{thm}

The condition of the existence of the pairing in (i) is satisfied  if $X$ is the reduction of $V$ as in Introduction, by the hard Lefschetz on $V$.
This assumption (i) is used in the proof of Lemma \ref{l3} below. 

\begin{lem}\label{l1} Let $v\in H^m(X)_{\ell}(r)^G$ with $m\geq 2r$. Then $v$ is in the image of $N^{m-2r}$.

\end{lem}

\begin{pf} 
By the  log Tate conjecture \ref{Tate4}, we have a $G$-homomorphism $\Sym^{m-2r}H^1(E)_{\ell} \to H^m(X)_{\ell}(r)$ which sends $e_1^{m-2r}$ to $v$ (see the last remark in \ref{Tatec}). Since $e_1^{m-2r}$ belongs to the image of $N^{m-2r}: \Sym^{m-2r}H^1(E)_{\ell}(m-2r)\to \Sym^{m-2r}H^1(E)_{\ell}$, $v$ belongs to the image of $N^{m-2r}: H^m(X)_{\ell}(m-r)\to H^m(X)_{\ell}(r)$. 
\end{pf}

\begin{para}\label{l0} To prove Theorem \ref{t:WM}, by \cite{IKNU} Proposition 2.1.13 (cf.\ \cite{Sa}, \cite{Vi} Proposition 2.4.2.1, \cite{Yo}), we are reduced to the case where $X$ is  strictly semistable. 

Then we can use the Rapoport--Zink spectral sequence in \cite{Nc:dege}.

\end{para}

\begin{lem}\label{l2}  Let $v\in H^m(X)_{\ell}^{N=0}$ and assume that $v$ is of weight $w\leq m$. Then $v$ is in the image of $N^{m-w}$.

\end{lem}

\begin{pf} We may assume $v\neq 0$ and $v$ is an eigenvector of Frobenius with eigenvalue $\alpha$. 
 Since $N(v)=0$, in the Rapoport--Zink spectral sequence (\cite{Nc:dege} Proposition 1.8.3), $v$ appears in $H^{2t+w}(Z)_{\ell}(t)$ for some projective smooth scheme $Z$  over $k$ (in the non-log sense) and for some integer $t$.  By Poincar\'e duality, the Frobenius eigenvalue $\alpha^{-1}$ appears in $H^{2t-w}(Z)_{\ell}(t)$ for some $Z$ and $t$ with the eigenvector $v'$.  Then $v\otimes v'\in 
H^{m+2t-w}(X\times Z)_{\ell}(t)^G$. Note that $m+2t-w\geq 2t$. By Lemma \ref{l1}, $v\otimes v'$ belongs to the image of $N^{m-w}$, that is, to the image of $N^{m-w}\otimes 1$ on $H^m(X)_{\ell} \otimes H^{2t-w}(Z)_{\ell}(t)$. Hence $v$ belongs to the image of $N^{m-w}$. 
\end{pf}

\begin{para}\label{l3} We complete the proof of Theorem \ref{t:WM}.

  By the existence of the pairing in  (i), the determinantial weight (the average weight) of $H^m(X)_{\ell}$ is $m$. 
  We prove that this together with Lemma \ref{l2} proves Theorem \ref{t:WM}.
  Fix a lifting of the Frobenius to $\pi_1^{\log}(s)$. 
  By the condition (ii), we have a decomposition $H^m(X)_{\ell}^{N=0}=\bigoplus_w V_w$, where $V_w$ is the part of weight $w$. 
  Take a base $(v_{w,i})_i$ of each $V_w$ for $w \le m$. 
  By Lemma \ref{l2}, for each $w \le m$ and $i$, there is a $v'_{w,i}\in H^m(X)_{\ell}$ such that $v_{w,i}=N^{m-w}(v'_{w,i})$.
  We may assume that this $v'_{w,i}$ is of  weight $w+2(m-w)= 2m-w$. 
  For $w \le m$ and $0 \le j \le m-w$, let $V_{w,+2j}$ be the subspace of $H^m(X)_{\ell}$ generated by $(N^{m-w-j}(v'_{w,i}))_i$ of  weight $w+2j$.  
  Thus $V_{w,+0}=V_w$. 
  Let $V=\bigoplus_{w \le m, 0 \le j \le m-w} V_{w,+2j}$. 
  It is sufficient to prove $H^m(X)_{\ell}=V$. 
  The average weight of $V$ is $m$ by construction. 
  On the other hand, we see  that the  weights of $H^m(X)/V$ are $>m$. 
  If $H^m(X) \neq V$, this contradicts that the average weight of $H^m(X)$ is $m$. 
\end{para}

\section{Log $K_0$ and the finer log Tate conjecture}\label{s:logK}

\begin{para}

In the classical theory, for a projective smooth scheme $X$ over a finite field, it is conjectured that $\Q_{\ell}\otimes \gr^rK_0(X)  \overset{\cong}\to H^{2r}(X)_{\ell}(r)^G$. 
  That is, not only the surjectivity, the isomorphism is conjectured. 
For the log case, $K_{0, \lim}$ is too big to have the isomorphism even  in the case where $k$ is finite. In fact, if $X$ is a curve, for each log modification  along a singular point of $X$, $\dim_{\Q}\Q \otimes$ Pic $=\Q\otimes \gr^1K_0$ increases by one, and hence $\Q\otimes \gr^1 K_{0,\lim}$ usually becomes of infinite dimension. To improve this situation, the following log $K_0$ may be nice. 

\end{para}

\begin{para}

For an element $a$ of the inductive limit of $H^0(X', M^{\gp}_{X'}/\cO^\times_{X'})$ for log modifications $X'$ of $X$, let $l_a\in K_{0,\lim}(X)$ be the class of the line bundle defined by the image of $a$ under the connecting map $H^0(X', M^{\gp}_{X'}/\cO^\times_{X'})\to H^1(X', \cO^\times_{X'})$. Let 
$$K_0^{\log}(X)$$ 
be the quotient of the ring $K_{0,\lim}(X)$ by the ideal generated by $l_a-1$ for all $a$. Then the Chern class maps to the $\ell$-adic cohomology factor through $K_0^{\log}(X)$. This is because the Chern character of $l_a$ in $\bigoplus_r \; H^{2r}(X)_{\ell}(r)$ is $1$ by the following lemma, and since the Chern character is a ring homomorphism, it kills the ideal generated by all $l_a-1$. Hence the Chern class maps kill this ideal.

\end{para}

\begin{lem} 

The Chern class of $l_a$ in $H^2(X)_{\ell}(1)$ vanishes. 

\end{lem}

\begin{pf}

This is because the Chern class map on $\text{Pic}$ factors through $H^1(M^{\gp})$ as is seen by the following commutative diagram.
$$\begin{matrix} 0&\to &\Z/\ell^n\Z(1) & \to & \mathbb{G}_m & \overset{\ell^n}\to & \mathbb{G}_m & \to & 0\\
&& \Vert && \cap &&\cap && \\
0&\to &\Z/\ell^n\Z(1) & \to & \mathbb{G}_{m,\log} & \overset{\ell^n}\to & \mathbb{G}_{m,\log} & \to & 0.\end{matrix}$$
\end{pf}

\begin{conj}\label{conjiso1} In the case where $k$ is a finite field, we have
$$\Q_{\ell} \otimes \gr^rK^{\log}_0(X) \overset{\cong}\to H^{2r}(X)_{\ell}(r)^G.$$

\end{conj}

\begin{para}\label{Kmr}  To consider the general $m\geq 2r$, we introduce a subgroup $K(X, m, r)$ of $\Q\otimes \gr^{m-r}K^{\log}_0(X\times E^{m-2r})$.

For an integer $a\geq 1$, the multiplication by $a: E^{(q)}\to E^{(q)}$ induces a morphism (also denoted by $a$) $E'=E^{(q, an)} \to E=E^{(q,n)}$.  For $1 \le i \le t$, 
by the pullback via the morphism $X\times (E')^t \to X\times E^t$ induced by the evident morphisms on all components except $a$ on the $i$-th $E'$, 
we have a homomorphism $a^*_i:   K_0^{\log}(X\times E^t) \to K_0^{\log}(X \times E^t)$.

Define $K(X,m,r)$ to be the part of $\Q \otimes \gr^{m-r}K^{\log}_0(X\times E^{m-2r})$ consisting of all elements $x$ satisfying the following conditions (i) and (ii).

\medskip

(i) $a^*_ix= ax$ for $1\leq i\leq m-2r$ and for any $a \ge1$.

(ii) $x$ is invariant under the action of the Symmetric group $S_{m-2r}$. 

\medskip

We define $K(X, m,r)=0$ in the case $m<2r$. 

We have $K(X, 2r,r)=\Q \otimes \gr^rK_0^{\log}(X)$. 
\end{para}

\begin{conj}\label{conjiso2}  In the case where $k$ is a finite field, we have
$$\Q_{\ell} \otimes K(X, m, r) \overset{\cong}\to H^m(X)_{\ell}(r)^G.$$
\end{conj}

\begin{para} We have the following relations of the conjectures.
Conjecture \ref{conjiso1} is a finer version of Conjecture \ref{Tate2}, and  
Conjecture \ref{conjiso2} is a finer version of Conjecture \ref{Tate4}. 

The idea of Conjecture \ref{conjiso2} comes from Conjecture \ref{conjiso1} as follows.
By Conjecture \ref{conjiso1},we expect
$$\Q_{\ell} \otimes \gr^{m-r}K^{\log}_0(X\times E^{m-2r})\overset{\cong}\to H^{2m-2r}(X\times E^{m-2r})_{\ell}(m-r)^G.$$
Consider elements of the right-hand-side satisfying the conditions (i) and (ii). By the condition (i) and by the fact that the pullbacks  $a^*: H^t(E)_{\ell} \to H^t(E)_{\ell}$ for integers $a\geq 1$ and $t\geq 0$  coincide with the multiplication by $a^t$, we get the K\"unneth component
$$H^m(X)_{\ell}\otimes H^1(E)_{\ell}^{\otimes (m-2r)}(m-r)^G.$$
By the condition (ii), we get 
$$\Hom_G(\Sym^{m-2r}H^1(E)_{\ell}, H^m(X)_{\ell}(r)).$$
If we admit the weight-monodromy conjecture, the last thing should be equal to $H^m(X)_{\ell}(r)^G$. 

\end{para}

\begin{para} Let $V$ be as in Introduction with strictly semistable reduction $X$, assume that $k$ is a finite field, and let $$L(H^m(V), s)= \text{det}(1-\varphi \cdot\sharp(k)^{-s}\;;\; H^m(V_{\bar K}, \Q_{\ell})^{N=0})^{-1}$$ be the Euler factor of $L$, where $\varphi$ is the geometric Frobenius. In \cite{CC}, a conjecture of S.\ Bloch on the order of the pole of $L(H^m(V), s)$ at $s=r\in \Z$ is discussed. He used his higher Chow group in the conjecture. We can give a modified version of his conjecture: 
\end{para}

\begin{conj}

The order of the pole of $L(H^m(V), s)$ at $s=r\in \Z$ is equal to $\dim_{\Q}(K(X, m, r))$.
\end{conj}

In fact, this order of the pole is equal to the multiplicity of the eigenvalue $1$ in the action of $\varphi$ on $H^m(V_{\bar K}, \Q_{\ell}(r))^{N=0}= H^m(X)_{\ell}(r)^{N=0}$. Hence if we assume the semisimplicity of the action of $\varphi$ as is widely believed, 
the log Tate conjecture \ref{conjiso2} implies this conjecture.

\begin{para}\label{int} The authors expect that an advantage of  our method is that by working with $K_0$, we can apply the usual method of algebraic cycles, for example, the intersection theory. The rationality and the positivity of the intersection number may be meaningful. 

We consider the following new intersection theory.

Assume $m,m', r, r'\in \Z$, $m, m'\geq 0$, $m+m'=2d$ ($d=\dim(X)$), $m-2r=m'-2r'\geq 0$. Then we have the pairing
$$K(X, m, r)\times K(X, m',r') \to \Q$$
 obtained as (letting $t=m-2r=m'-2r'$) 
$$\gr^{m-r}K_{0,\lim}(X \times E^t) \times \gr^{m'-r'}K_{0,\lim}(X \times E^t) \to \gr^{(m+m')-(r+r')}K_{0,\lim}(X \times E^t)$$
$$= \gr^{d+t}K_{0,\lim}(X \times E^t)\to \Q,$$
where the last arrow is the left vertical arrow of the diagram in 
 \cite{IKNU} Proposition 2.4.9 which we use by taking $X \times E^t$ and $s$ here as $X$ and $Y$ there, respectively.

\end{para}

\begin{conj} This pairing 
$K(X, m,r)\times K(X, m',r') \to \Q$
is a perfect pairing of finite dimensional $\Q$-vector spaces.

\end{conj}

\begin{para} This pairing has the following property. If $a\in K(X, m,r)$, $b\in K(X, m', r')$, and if $\alpha\in H^m(X)_{\ell}(r)^G$ and $\beta\in H^{m'}(X)_{\ell}(r')^G$ are the images of $a$ and $b$, respectively, the pairing $(a,b)\in \Q$ coincides with $\tilde \alpha\cup \beta$, where $\tilde \alpha$ is an element of $H^m(X)_{\ell}(m-r)$ such that $N_q^{m-2r}(\tilde \alpha)=\alpha$ and $\cup$ is the cup product $H^m(X)_{\ell}(m-r) \times H^{m'}(X)_{\ell}(r') \to H^{2d}(X)(m-r+r') =H^{2d}(X)_{\ell}(d) \to \Q_{\ell}$. Here $q$ in $N_q$ is the $q$ which we used for  $E=E^{(q,n)}$. 

In \ref{HHQ}, we consider an example of this. 
\end{para}

 \section{The monodromy cycle}\label{s:mod}

\begin{para}\label{Gm} Assume that $k$ is a finite field.
If we use $X\times X$ as $X$ in Conjecture \ref{conjiso2} and use the K\"unneth decomposition and Poincar\'e duality on the log \'etale cohomology, we would obtain an isomorphism
$$\Q_{\ell}\otimes K(X\times X, 2d, d-1) \overset{\cong}\to 
   {\textstyle\bigoplus}_m \; \Hom_G(H^m(X)_{\ell}, H^m(X)_{\ell}(-1)),$$
where $d=\dim X$.  

This suggests that for a section $q$ of the log structure of $s$ which does not belong to $k^\times$, we have a unique element of $K(X \times X, 2d, d-1)$  which induces the monodromy operators $N_q: H^m(X)_{\ell}\to H^m(X)_{\ell}(-1)$ for all $m$ and all $\ell\neq p$ and also the analogous monodromy operator of the log crystalline cohomology.
  We call this expected element the {\it monodromy cycle}. 

\end{para}

\begin{conj}
  A monodromy cycle  exists. 
\end{conj}

\begin{para}\label{sl2} A consequence of the weight-monodromy conjecture is that (choosing a lifting of the geometric Frobenius from $\Gal(\bar k/k)$ to $\pi_1^{\log}(s)$ and choosing a base of $\Q_{\ell}(1)$) we have an action of the Lie algebra $\sl(2)$ on $H^m(X)_{\ell}$ characterized by the following properties. The monodromy operator $N$ on $H^m(X)_{\ell}$ is the action of the matrix $\begin{pmatrix} 0&1\\0&0\end{pmatrix}\in \sl(2)$, and the matrix
$\begin{pmatrix} -1 & 0 \\ 0 & 1 \end{pmatrix}\in \sl(2)$ acts on the part of $H^m(X)_{\ell}$ of  Frobenius weight $w$ as the multiplication by $w-m$. 

In the following \ref{mon2}, we give a geometric understanding of this  action of $\sl(2)$  by using the monodromy cycle. 

In \ref{try}, we will try to understand the monodromy cycle and this action of $\sl(2)$ by using the theory of log motives \cite{IKNU}. 

\end{para}

\begin{para}\label{mon2}

 The map
$$\gr^{d+1}K_0^{\log}(X \times X \times E \times E) \to H^{2(d+1)}(X \times X\times E \times E)_{\ell}(d+1)$$
induces, by the condition (i) in \ref{Kmr}, a homomorphism
$$K(X \times X, 2d, d-1) \to (H^{2d}(X \times X)_{\ell}(d) \otimes H^1(E)_{\ell} \otimes H^1(E)_{\ell}(1))^G$$ 
$$= \Hom_G(\End(H^1(E)_{\ell}), {\textstyle\bigoplus}_m  \;\End(H^m(X)_{\ell})),$$ 
and the condition (ii) in \ref{Kmr} tells that the image of this map is contained in the space of homomorphisms $\End(H^1(E)_{\ell})\to  {\textstyle\bigoplus}_m  \;
\End(H^m(X)_{\ell})$ which kill the scaler matrices $\Q_{\ell}$ in $\End(H^1(E)_{\ell})$. Note that 
$\End(H^1(E)_{\ell})/\Q_{\ell}$ is identified with the part $\sl(H^1(E)_{\ell})$ of $\End(H^1(E)_{\ell})$ of trace $0$.
Thus we have
$$K(X \times X, 2d, d-1) \to \Hom_G(\sl(H^1(E)_{\ell}), {\textstyle\bigoplus}_m  \;\End(H^m(X)_{\ell}))$$
If $[N]\in K(X \times X, 2d, d-1)$ is a monodromy cycle, the induced homomorphism $$\sl(H^1(E)_{\ell})\to {\textstyle\bigoplus}_m  \;\End(H^m(X)_{\ell})$$ 
is the action of $\sl(2)$ mentioned in \ref{sl2}. This homomorphism  sends the monodromy operator $N_q:H^1(E)_{\ell} \to H^1(E)_{\ell}(-1)$ of $E$ to the monodromy operators $N_q: H^m(X)_{\ell} \to H^m(X)_{\ell}(-1)$  of $X$. (Here this $q$ need not be equal to the $q$ which we used to define the log Tate curve $E$.)

\end{para}

\begin{rem}
In \cite{logmot} A.17, we conjectured the existence of a monodromy cycle as an element of $\gr^{d-1}KH_{-2,\lim}(X \times X)$ not assuming the finiteness of $k$. But this element does not  give an action of $\sl(2)$.  The action of $\sl(2)$ appears only by the power of the log Tate  curve. 
\end{rem}

\begin{para} 
In \cite{IKNU} Section 5.1, the category of log motives $LM_R(S)$ (the logarithmic version of the category of Grothendieck motives) over an fs log scheme $S$  was defined fixing a non-empty set $R$ of prime numbers which are invertible on $S$, by using the homological equivalence for $\ell$-adic realizations for $\ell\in R$.
We consider the case where $S$ is our log point $s$ and $k$ is a finite field, and $R$ is the set of all prime numbers which are not the characteristic of $k$. We denote $LM_R(s)$ by $\frak M=\frak M(s)$. 

Objects of $\frak M$ are 
direct summands of the objects $h(X)(r)$ which are associated to objects $X$ of $\frak S$ and integers $r$ in $\bZ$.
 For objects $X$ and $X'$ of $\frak S$ and for $r,r'\in \Z$, we have a surjective homomorphism
$$\Q\otimes \gr^{\dim(X)+r'-r}K_0^{\log}(X\times X')\to \Hom_{\frak M}(h(X)(r), h(X')(r')),$$
and this is an isomorphism if the finer $\ell$-adic log Tate conjecture \ref{conjiso1} is true for some $\ell\in R$.

\end{para}

\begin{para}\label{try1} Let
$$\cH om: \frak M\times \frak M \to \frak M, \quad \cE nd: \frak M \to \frak M$$
be internal hom and the internal end, respectively. Using the dual log motive $M^*$ of $M$, they are written as $(M, M') \mapsto M^*\otimes M'$ and $M\mapsto M^*\otimes M$, respectively. 
Let 
$\frak A=\cE nd (H^1(E))=H^1(E)^*\otimes H^1(E)$. Then $\frak A$  is a ring object of $\frak M$. It is hence a Lie algebra object with the commutator Lie bracket.  Let $\sl(H^1(E))$ be the  kernel of the canonical morphism $H^1(E)^*\otimes H^1(E) \to 1$, where $1$ is the unit object of $\frak M$. Then $\sl(H^1(E))$ is a Lie algebra object of $\frak M$.

We try to understand the monodromy cycle of a log motive $M\in \frak M$  as a canonical homomorphism of Lie algebra objects $\sl(H^1(E)) \to \cE nd (M)$, that is, a canonical action of $\sl(H^1(E))$ on $M$ in $\frak M$. 

\end{para}

\begin{para}\label{try2} Consider the full subcategory $\frak M_{\sat}=\frak M_{\sat}(s)$ of $\frak M$ consisting of all objects which are isomorphic to a direct summand of $h(X)(r)$ for some $X\in \frak S$ which is saturated over $s$ and for some $r\in \Z$. 

Note that $X\in \frak S$ is saturated over $s$ if and only if $X$ is reduced as a scheme (\cite{Tsu} Theorem II.4.2). 

If $X$ and $Y$ are saturated over $s$, the fiber product $X\times Y$ over $s$ is saturated. Hence in $\frak M$, $\frak M_{\sat}$ is stable under tensor products, duals, and direct summands.

If $X\in \frak S$ is semistable, then  $X$ is saturated. Though we feel that the properties semistable and saturated are similar good properties,
we use the latter to define a subcategory in the above because the property semistable is not stable under the fiber products over $s$. 

\end{para}

\begin{prop}\label{satunip} 
Let $M$ be an object of $\frak M_{\sat}$. Then the action of the inertia subgroup $\pi_1^{\log}(\bar s)$ of $\pi_1^{\log}(s)$ on the $\ell$-adic realization $M_{\ell}$ of $M$ is unipotent.

\end{prop}

This follows from the next proposition. 

\begin{prop}\label{unip}

 Let $f: Y\to s$ be a morphism from an fs log scheme to the standard log point and $\ell$ a prime invertible on $s$.

  $(1)$  Assume that the cokernel of
  $$\Z=f^{-1}(M^{\gp}_s/\cO^\times_s)\to M_Y^{\gp}/\cO^\times_Y$$
  is torsion-free. Then for every $m\geq 0$ and $n\geq 0$, the action of every element $g$ of $\pi_1^{\log}(\bar s)$ on $H^m(Y_{\bar s(\log)}, \Z/\ell^n\Z)=R^mf_{\loget,*}(\Z/\ell^n\Z)_{\bar s(\log)}$ satisfies $(g-1)^{m+1}=0$. 
  
  $(2)$ The assumption on $f$ in $(1)$ is satisfied if $f$ is saturated.
  
  \end{prop}

\begin{pf} (2) is \cite{IKN} Proposition (A.4.1).

We prove (1). Let $\overset{\circ}Y$ be the underlying scheme of $Y$ endowed with the \'etale topology.
  By the Leray spectral sequence 
$$E_2^{i,j}=H^i_{\et}(\overset{\circ}Y,\Psi^j)\Rightarrow H^{i+j}_{\loget}(Y_{\overline s(\log)},\Z/\ell^n\Z),$$
it is enough to show that the action of $\pi_1^{\log}(\overline s)$ on the nearby cycle 
$$\Psi^j:=R^j(Y_{\overline s(\log)} \overset \pi \to Y \overset \varepsilon \to \overset{\circ}Y)_*(\Z/\ell^n\Z)
=R^j\varepsilon_*\pi_*\Z/\ell^n\Z$$
is trivial. 
  This is a local problem and can be checked stalkwise as follows.
  Let $y \in Y$. 
  Then the argument in the proof of \cite{Nc:dege} Lemma 1.8.2, which treats the semistable case, works if the homomorphism 
$f\colon \pi_1^{\log}(\bar y) \to \pi_1^{\log}(\bar s)$ of profinite groups has a section. For an fs log point $x$, $\pi_1^{\log}(\bar x)=\Hom(M^{\gp}_{\bar x}/\cO^\times_{\bar x}, \prod_{\ell'} \Z_{\ell'}(1))$, where $\ell'$ ranges over all prime numbers which are invertible on $x$. Hence the existence of the section follows from the fact that
   the cokernel of the injection $M_{\bar s}^{\gp}/\cO^{\times}_{\bar s} \to M_{\bar y}^{\gp}/\cO^{\times}_{\bar y}$ is torsion-free.
         \end{pf}

\begin{para}
For a log Tate curve $E$, the log motive $H^1(E)$ belongs to $\frak M_{\sat}$ because it is $H^1(E^{(q, n)})$ for $q$ and $n\geq 2$ such that the image of $q$ in $\N=M_s/\cO^\times_s$ is a multiple of $n$ and $E^{(q,n)}$ is saturated for such $(q,n)$. Hence for any $r\geq 0$ and any classical Grothendieck motive $C$ over $k$, $\Sym^rH^1(E) \otimes C$ belongs to $\frak M_{\sat}$.

Now Proposition \ref{satunip},  \cite{KNU} Proposition 1.10 and the log Tate conjecture \ref{conjiso2} suggest the following.
\end{para}

\begin{conj}\label{motconj} Let $M$ be an object of $\frak M_{\sat}$. Then 
 $M$ is isomorphic to $${\textstyle\bigoplus}_{r\geq 0}\; \Sym^r(H^1(E)) \otimes C_r$$ for some classical Grothendieck motives $C_r$ over $k$ which are $0$ for almost all $r$.

\end{conj}

\begin{para}\label{try}

For $M$ as in Conjecture \ref{motconj}, the desired action of $\sl(H^1(E))$ on $M$ is given as the tensor product of the natural action on $\Sym^r(H^1(E))$ and the trivial action on $C_r$. This action is 
well-defined (we assume Conjecture \ref{conjiso1} here). In fact, $C_r$ is canonical because it represents the functor $D\mapsto \Hom(\Sym^r(H^1(E)) \otimes D, M)$ on the category of classical Grothendieck motives over $k$, and the morphism $\Sym^r(H^1(E)) \otimes C_r\to M$ comes from this property of $C_r$, and hence the object $C_r$ and the isomorphism in Conjecture \ref{motconj} are canonical.

\end{para}

\section{Examples}\label{s:Ex}

Assume that $k$ is a finite field.
The next theorem proves a special case of Conjecture \ref{conjiso2}. 

\begin{thm}\label{prop1}

Assume that $X$ is a strictly semistable curve. 
Then the map
$\Q_{\ell}\otimes K(X, 1, 0)\to H^1(X)_{\ell}^G$ is surjective. 

\end{thm}

\begin{para}\label{pf1} We prove Theorem \ref{prop1} till \ref{pf4}.
 In \ref{pf4}, we use the additivity of the Albanese map, which will be seen in the next section by using theta functions. 

Let $I$ be the set of all irreducible components of $X$ and let $J$ be the set of all singular points of $X$. By using the base change to a finite extension of $k$, we may assume that all points in $J$ are $k$-rational. 
   We may assume that $X$ is connected.

Define the free abelian group $H^1$ by the exact sequence
$$0 \to \Z \to {\textstyle\bigoplus}_{i\in I} \;\Z \to {\textstyle\bigoplus}_{j\in J} \; \Z \to H^1\to 0.$$
This is identified with the exact sequence of  cohomology groups for $*=\zar$ or \'et
$$0\to H^0_*(X, \Z)\to H^0_*(U, \Z)\to {\textstyle\bigoplus}_{j\in J}\;  H^1_{*,j}(X, \Z)\to H^1_*(X, \Z) \to H^1_*(U, \Z),$$
where $U$ is the non-singular part of $X$ and $H^1_{*,j}(X,\cdot)$ is the cohomology with support. That is, $$H^1= H^1_{\zar}(X, \Z)=H^1_{\et}(X, \Z).$$
Take the $\Z$-dual of the above exact sequence 
$$0 \to H_1\to {\textstyle\bigoplus}_{j\in J}\; \Z \to {\textstyle\bigoplus}_{i\in I}\; \Z \to \Z\to 0,$$
where $H_1$ is the $\Z$-dual of $H^1$. 

Then the theory of the Rapoport--Zink spectral sequence (\cite{Nc:dege}) gives the following.
$H^1(X)_{\ell}$ has an increasing filtration $W$ such that $W_{-1}=0$, $W_2=H^1(X)_{\ell}$, $W_0=\Q_{\ell}\otimes H^1$, $W_1/W_0={\textstyle\bigoplus}_{i\in I} \;H^1(X_i)_{\ell}$, where $X_i$ denotes the irreducible component $i\in I$, $W_2/W_1= \Q_{\ell}(-1)\otimes H_1$, and the monodromy operator $N_q$ for a generator $q$ of the log structure of $s$ coincides with the composition 
$$H^1(X)_{\ell}\to \Q_{\ell}(-1) \otimes H_1 \to \Q_{\ell}(-1) \otimes H^1 \to H^1(X)_{\ell}(-1),$$
in which the middle arrow is induced by the composition of the canonical maps
$$H_1 \to {\textstyle\bigoplus}_{j\in J}\; \Z \to H^1.$$
Though we are working over a log point $s$, not over a discrete valuation ring, by the same argument as in the case of curves over a discrete valuation ring, we have that the map $\Q\otimes H_1 \to \Q \otimes H^1$ is an isomorphism, that is, the weight-monodromy conjecture holds for this $H^1(X)_{\ell}$.

\end{para}

\begin{lem}\label{pf2} Let $q$ be a section of the log structure of $s$ which does not belong to $k^\times$. Then the image of  $H^1_{\zar}(X, \Z)\to H^1_{\zar}(X, M^{\gp}_X)$ induced by the homomorphism $\Z\to M^{\gp}_X\;;\;1\mapsto q$ is finite.

\end{lem}

\begin{pf}
We use the result in Kajiwara \cite{Kaj}. Working on the Zariski site of $X$, consider the quotient sheaf 
 $\cF= M_X^{\gp}/(\cO_X^\times \oplus q^{\Z})$. Then $\cF\cong {\textstyle\bigoplus}_{j\in J}\; (a_j)_*\Z$, where $a_j$ is the inclusion map $j\to X$. From the exact sequence $0\to \cO_X^\times \oplus q^{\Z} \to M_X^{\gp} \to \cF \to 0$, 
 we obtain an exact sequence of cohomology
 $${\textstyle\bigoplus}_{j\in J} \; \Z \to \text{Pic}(X) \oplus H^1 \to H^1_{\zar}(X, M_X^{\gp})\to 0.$$
By \cite{Kaj} Theorem 2.19, this exact sequence induces an exact sequence of subgroups 
$$0\to H_1 \to \text{Pic}(X)_0 \oplus H^1 \to H^1_{\zar}(X, M_X^{\gp})_0\to 0.$$
  Here $\text{Pic}(X)_0$ is the kernel of $\text{Pic}(X) \to {\textstyle\bigoplus}_{i\in I} \; \text{Pic}(X_i) \to {\textstyle\bigoplus}_{i\in I} \; \Z$, where the second arrow is given by the degree maps, and $H^1_{\zar}(X, M^{\gp}_X)_0$   is the kernel of a homomorphism $H^1_{\zar}(X, M^{\gp}_X) \to \Z$ which is called the degree map in \cite{Kaj}. Consider the exact sequence of sheaves $0\to \cO_X^\times \to {\textstyle\bigoplus}_{i\in I} \; (p_i)_*(\cO_{X_i}^\times)\to \cG\to 0$, where $p_i$ is the inclusion morphism $X_i\to X$, and the stalk of $\cG$ at $x\in X$ is $k^\times$ if $x$ is a singular point and is $0$ otherwise. By the finiteness of $k$, we have that
the kernel of $\text{Pic}(X) \to {\textstyle\bigoplus}_{i\in I} \; \text{Pic}(X_i)$ is finite. Since the kernel of the degree map $\text{Pic}(X_i)\to \Z$ is finite by the finiteness of $k$, we have that $\text{Pic}(X)_0$ is finite. Hence by the above \cite{Kaj} Theorem 2.19, we have an exact sequence $0\to \Q \otimes H_1 \to \Q \otimes H^1 \to \Q \otimes H^1_{\zar}(X, M^{\gp}_X)_0 \to 0$. Since $\Q \otimes H_1 \to \Q \otimes H^1$ is an isomorphism, the map $\Q \otimes H^1 \to \Q \otimes H^1_{\zar}(X, M^{\gp}_X)$  is the zero map.\end{pf}

\begin{lem}\label{pf3}
Consider the exact sequence $0 \to \Z \overset{q}\to \mathbb{G}_{m,\log}^{(q)}\to E^{(q)}\to 0$ on $X$. 
The connecting map $$\Phi: \Q\otimes_{\Z} E^{(q)}(X) \to \Q \otimes H^1$$  is surjective.

\end{lem}

\begin{pf}
By Lemma \ref{pf2}, 
we have the surjection $\Q\otimes H^0(X, \mathbb{G}_{m,\log}/q^{\Z})\to \Q \otimes H^1$. Since $X$ is vertical, $H^0(X, \mathbb{G}_{m,\log}/q^{\Z})= H^0(X,  \mathbb{G}_{m,\log}^{(q)}/q^{\Z})= E^{(q)}(X)$.
\end{pf}

\begin{para}\label{pf4} 

By \cite{IKNU} Section 3.1, for $a\in E^{(q)}(X)$, we have an element $C(a)$ of $\gr^1K_{0,\lim}(X \times E)$ whose image under 
$\gr^1K_{0,\lim}(X \times E) \to H^1(X)_{\ell}^G$ coincides with $\Phi(a)$. This element is obtained as follows (see the proof of Proposition 3.1.4 in \cite{IKNU}). 

 Let $E=E^{(q,n)}$ and let $Y$ be the log modification of $E\times E$ defined as $Y=\tilde Y/(q^{\Z}\times q^{\Z})$, where 
 $\tilde Y$ is the log modification of  $\mathbb{G}_{m,\log}^{(q,n)} \times \mathbb{G}_{m,\log}^{(q,n)}$ which represents the following functor. For an fs log scheme $T$ over $s$, $\tilde Y(T)$ is the set of all $(t_1, t_2)\in (\mathbb{G}_{m,\log}^{(q,n)} \times \mathbb{G}_{m,\log}^{(q,n)})(T)$ such that for each $r\in \Z$, we have either $q^rt_1|t_2$ or $t_2|q^rt_1$ locally on $T$. 
 Then the
diagonal $\Delta=E\to E\times E$ factors through a strict closed immersion $\Delta\to Y$. The ideal of $\cO_Y$ which defines $\Delta$ is an invertible and hence is a line bundle on $Y$.
 Let $Z$ be a log modification of $X\times E$ such that the composition $Z\to X\times E \overset{(a,\text{id.})}\longrightarrow E^{(q)} \times E$ factors as $Z\overset{a'}\to Y \to E \times E \to E^{(q)}\times E$.  Let $C(a)\in \gr^1K_{0,\lim}(X\times E)$ be the pullback of the class of this line bundle in $\gr^1K_0(Y)$ under $a'$.

Let $\Psi(a)$ be the image of $C(a)-C(0)$ in $\gr^1K_0^{\log}(X \times E)$. Since the image of $C(a)$ in $H^1(X)_{\ell}^G$ is $\Phi(a)$ and the image of $C(0)$ in $H^1(X)_{\ell}^G$ is $\Phi(0)=0$, we have:

\medskip

(1) The following diagram is commutative. 
$$\begin{matrix}  E^{(q)}(X) &\overset{\Psi}\to & \gr^1K_0^{\log}(X \times E)\\
\overset{\Phi}{}\downarrow{}{}{} && \downarrow\\
\ \ H^1 &\to& H^1(X)_{\ell}^G\end{matrix}$$

We will prove the following 
(2)  in \ref{theta4} below.

\medskip

(2) The map $\Psi: \text{Mor}(X, E^{(q)})\to\gr^1K_0^{\log}(X \times E)$ is a homomorphism.

\medskip

\noindent (Though this (2) can be deduced from the theory of the log Jacobian variety of $X$ and its self-duality described in \cite{IKNU} Section 6.2, we give a proof of (2) in Section \ref{s:theta} below  using the theta function of a log Tate curve over our log point $s$ because we think the log Tate curve is important in this paper and its theta function may have a value and because the full details of \cite{IKNU} Section 6.2 are not yet published.)
By (2), we have

\medskip

(3) The image of $\Psi$ in $\Q\otimes \gr^1K_0^{\log}(X \times E)$ is contained in $K(X, 1,0)$. 

\medskip

 By Lemma \ref{pf3} and the above (1) and  (3), we have Theorem \ref{prop1}.

\end{para}

\begin{para}\label{HHQ} Let $X$ be as in Theorem \ref{prop1} and assume that all singular points are $k$-rational. Take the Tate elliptic curve $E^{(q,n)}$, where $q$ is a generator of the log structure of $s$. Then we have a commutative diagram 
$$\begin{matrix}  K(X, 1, 0) \times K(X, 1,0) &\to & \Q\\
\downarrow && \Vert
\\ \Q\otimes H^1 \times \Q \otimes H^1 & \to & \phantom{.}\Q.\end{matrix}$$
Here the upper arrow is as in \ref{int} and the lower arrow is $(\Q\otimes H^1)\otimes (\Q\otimes H^1)  \cong (\Q\otimes H_1)\otimes (\Q\otimes H^1) \to \Q$ in which the first arrow is by the isomorphism $\Q\otimes H_1 \overset{\cong}\to \Q\otimes H^1$ in \ref{pf1}.

\end{para}

The next theorem is about the existence of a monodromy cycle in a special case.

\begin{thm}\label{Thmon} Assume that $X$ is a strictly semistable curve.

$(1)$ The map $\Q_{\ell}\otimes K(X\times X, 2, 0)\to H^2(X \times X)_{\ell}^G={\textstyle\bigoplus}_m \Hom_G(H^m(X)_{\ell}, H^m(X)_{\ell}(-1))$ is surjective.

$(2)$ There is an element of $K(X\times X, 2, 0)$ which induces a homomorphism 
$$\sl(H^1(E)) \to\cE nd(h(X))$$
of Lie algebra objects in the category $\frak M$, 
which sends the monodromy operator $N_q: H^1(E)_{\ell}\to H^1(E)_{\ell}(-1)$ to the monodromy operator $N_q: H^m(X)_{\ell}\to H^m(X)_{\ell}(-1)$ for every $m$ and every $\ell$. 
\end{thm}

\begin{pf} We may assume that all singular points of $X$ are $k$-rational.

(1) is proved as follows. By \ref{pf1}, $H^2(X\times X)_{\ell}^G= H^1(X)_{\ell}^G \otimes H^1(X)_{\ell}^G$ and it is generated by the image of $H^1 \otimes H^1$. By Theorem \ref{prop1}, it is generated by the image of
$K(X, 1,0) \times K(X, 1, 0) \to K(X\times X, 2, 0) \to H^2(X \times X)_{\ell}^G$.

 (2) is proved as follows. 
 
 The surjective homomorphism $\Q\otimes E^{(q)}(X) \to \Q\otimes H^1$ and the homomorphism $\Psi: E^{(q)}(X) \to \gr^1K_0^{\log}(X\times E)$ give a homomorphism $H^1(E) \otimes H^1\to H^1(X)$ of log motives in $\frak M$, where $H^1(E) \otimes H^1$ means the log motive $H^1(E)^{\oplus r}$ if we fix a $\Z$-base $(e_i)_{1\leq i\leq r}$ of $H^1$. By Poincar\'e duality, this gives a homomorphism $H^1(X) \to H^1(E) \otimes H_1$ of log motives, and the composition 
 $H^1(E)\otimes H^1 \to H^1(X) \to H^1(E) \otimes H_1\to H^1(E)\otimes H^1$ is the identity morphism. Hence $H^1(E) \otimes H^1$ is regarded as a direct summand of the log motive $h(X)$. We define the Lie action of $\sl(H^1(E))$ on $h(X)$ by using the evident action on $H^1(E)$ on this direct summand, and using the trivial action on the other direct summand. This gives the monodromy operator as stated in (2). 
 
 This Lie action is induced by an element of $K(X\times X, 2, 0)$ obtained as follows. 
 The pairing $(\Q\otimes H^1)\times (\Q\otimes H^1)\to \Q$ in \ref{HHQ} is perfect. Take elements $a_1, \dots, a_n$ of $K(X, 1,0)$ such that if $\alpha_i$ denotes the image of $a_i$ in $\Q\otimes H^1$, then $(\alpha_i)_i$ is a base of $\Q\otimes H^1$. Let $b_1, \dots, b_n$ be elements of $K(X,1,0)$ such that if $\beta_i$ denotes the image of $b_i$ in $\Q\otimes H^1$, then $(\beta_i)_i$ is the dual base of $(\alpha_i)_i$ for the pairing $(\Q\otimes H^1)\times (\Q\otimes H^1) \to \Q$. The desired element of $K(X\times X, 2, 0)$ is obtained as  $\sum_{i=1}^n a_i\cdot b_i$, where $\cdot$ is the pairing $K(X, 1,0) \times K(X, 1,0) \to K(X\times X, 2, 0)$. 
 \end{pf}

\begin{para} In this paper, if we drop the assumption that $X$ is projective, many points  stop to work. 

For example, consider the reduction of a log Hopf surface. A Hopf surface over a $p$-adic local field $K$ is a rigid analytic space and not algebraic, but can have a reduction $X$ which is a proper log smooth scheme over $s$ of strictly semistable reduction, and which is not projective. 

 We have that $H^3(X)_{\ell}(2)^G=H^3(X)_{\ell}(2)=\Q_{\ell}$  is non-zero, so Conjecture \ref{Tate1} is not true for this $X$. 
 
 We have
$H^1(X)_{\ell}^G=H^1(X)_{\ell}=\Q_{\ell}$ but it does not come from a $G$-homomorphism $H^1(E)_{\ell}\to H^1(X)_{\ell}$.  So Conjecture \ref{Tate4} is not true for this $X$. 
\end{para}

 \section{Theta functions}\label{s:theta}

 Here we describe the theory of theta functions on a log Tate curve over a standard log point. This is the reduction to $k$ of the theory of theta functions on the Tate curve over a $p$-adic local field $K$ described, for example,  in \cite{Ro}. The latter theory is the $p$-adic analogue of the classical theory over $\C$ reviewed in \ref{theta1} below.
 
  This is used in \ref{pf4} in the previous section.
 
 \begin{para}\label{theta1} Let $q\in \C^\times$, $|q|<1$, and consider the elliptic curve $E=\C^\times/q^{\Z}$ over $\C$. We have the theta functions $\theta$ and $\theta_a$ for $a\in \C^\times$, which are meromorphic functions on $\C^\times$, defined as follows:
 $$\theta(t)= \prod_{n=0}^{\infty}\; (1-q^nt) \cdot \prod_{n=1}^{\infty}\; (1-q^nt^{-1}), \quad\quad  \theta_a(t)= \theta(t/a).$$
 Hence $\theta=\theta_1$. 
  We have
 
 \medskip
 
 (1) $\theta(qt)= -t^{-1}\theta(t), \quad \theta_a(qt)=-at^{-1}\theta_a(t).$

 (2) The divisor of $\theta_a$ on $\C^\times$ is the pullback of the divisor $(a)$ on $E$. 

\medskip

 From (1), we have:
 
 \medskip
 
 (3)  For $a,b\in \C^\times$, the function $f(t)=\theta_{ab}(t)\theta_1(t)\theta_a(t)^{-1}\theta_b(t)^{-1}$ satisfies $f(qt)=f(t)$ and hence is a meromorphic function on $E=\C^\times/q^{\Z}$. 
 
 By (2), the divisor of $f$ is the divisor $(ab)+(1) - (a)-(b)$ on $E$. This gives a proof of the fact that the divisor $(ab)+(1) -(a)-(b)$ on $E$ is a principal divisor.

 \end{para}
 
 \begin{para}\label{RQ} We now go to the log geometry.

 In general, for an fs log scheme $S$, let $\cR_S$ be the sheaf of rings over $\cO_S$ obtained from $\cO_S$ by inverting all non-zero-divisors in $\cO_S$. Let $\cQ_S$ be the log structure on $\cR_S$ associated to the pre-log structure $M_S \to \cR_S$. 
 
  The group sheaf $\cQ^{\gp}_S$ is generated by the subgroup sheaves $\cR_S^{\times}$ and $M^{\gp}_S$.
 We regard $\cQ_S^{\gp}/M_S^{\gp}$ as the log version of the sheaf $\cR_S^\times/\cO_S^\times$ of Cartier divisors. An element $D$ of $H^0(S, \cQ^{\gp}_S/M_S^{\gp})$ gives an $M^{\gp}_S$-torsor on $S$ defined to be the inverse image of $D$ in $\cQ^{\gp}_S$. It is like the $\cO_S^\times$-torsor associated to an element of  $H^0(S, \cR_S^\times/\cO^\times_S)$. We have the commutative diagram
 $$\begin{matrix}  H^0(S, \cR^\times_S/\cO^\times_S) & \to & H^1(S, \cO_S^\times)\\
 \downarrow && \downarrow\\
 H^0(S, \cQ_S^{\gp}/M^{\gp}_S)&\to& H^1(S, M_S^{\gp}).\end{matrix}$$ 
 
 \end{para}
 
 \begin{para}\label{theta2}  Now we work over a standard log point $s$.

 Let $E=E^{(q,n)}$ with $n\geq 2$ and let $\tilde E=\mathbb{G}_{m,\log}^{(q,n)}$ (\ref{Eqn}), so $E=\tilde E/q^{\Z}$. 
  
 Then there is a unique section $\theta=\theta(t)$  $(t \in \tilde E \subset \Gmlog^{(q)})$ of $\cQ^{\gp}_{\tilde E}$ satisfying the following conditions (i) and (ii).
 
 \medskip
 
 (i) $\theta(qt)= -t^{-1}\theta(t)$. 
 
 (ii) On the open part $\{t \in \Gmlog\,|\,1|t|q\}$ of $\tilde E$, $\theta=(1-\alpha(t))(1-\alpha(qt^{-1}))\in \cR_{\tilde E}^\times$, where $\alpha$ is the structural map $M \to \cO$ of the log structure. 
\medskip

The restriction of $\theta(t)$ to $\mathbb{G}_m\subset \tilde E$ is $1-t$.

 \end{para}
 
 \begin{para}\label{theta3}  Let $X$ be an object of $\frak S$ and let $a\in E^{(q)}(X)$. We will define the theta function $\theta_a$ as a section of $\cQ^{\gp}$ of a certain covering $\tilde E_{X,a}$ of a log modification $E_{X,a}$ of $X\times E$.
 
 First we consider the case $X=E$ and $a$ is the identity map $\delta: E\to E$. In this case, $E_{E,\delta}=Y$ and $\tilde E_{E,\delta}=\tilde Y$ in 
 \ref{pf4}. 
 The restriction of $\theta_{\delta}$ to $\mathbb{G}_m\times \mathbb{G}_m\subset \tilde Y$ is $1-t_1^{-1}t_2$. We define $\theta_{\delta}$ as the unique section of $\cQ^{\gp}$ of $\tilde Y$ characterized by the following properties (i) and (ii).
 
 \medskip
 
 (i)  $\theta_{\delta}(qt_1, qt_2)=\theta_{\delta}(t_1, t_2)$ and $\theta_{\delta}(t_1, qt_2)= -t_1t_2^{-1}\theta_{\delta}(t_1, t_2)$. 
 
 (ii) On the part $\{(t_1, t_2)\;|\; 1|t_1|q,\; 1|t_2|q\}$ of $\tilde Y$, $\theta_{\delta}(t_1, t_2)= (1-\alpha(t_1^{-1}t_2))(1-\alpha(qt_2^{-1}t_1))\in \cR^\times$ if $t_1|t_2$, and $\theta_{\delta}(t_1, t_2)= -t_1^{-1}t_2(1-\alpha(t_2^{-1}t_1))(1-\alpha(qt_2^{-1}t_1))\in \cQ^{\gp}$ if $t_2|t_1$.
 \medskip
 
 The restriction of $\theta_{\delta}$ to $\{1\}\times \tilde E$ is $\theta$. 
 
 \end{para}
 
 \begin{lem}\label{theta8} By the canonical map $\cR^\times/\cO^\times \to \cQ^{\gp}/M^{\gp}$ on $\tilde Y$, the Cartier divisor of the diagonal $E$ in $Y$ is sent to the class of $\theta_{\delta}$.

 \end{lem}
 
 \begin{pf} This is clear by the explicit form of $\theta_{\delta}$. 
 \end{pf}
 
 \begin{para}\label{theta7} Let $X$ be an object of $\frak S$ and let $a\in E^{(q)}(X)$. We define a log modification $E_{X,a}$ of  $X\times E$ to be the fiber product of $X\times E\overset{(a,1)}\to E^{(q)} \times E\leftarrow Y$. We define $\tilde E_{X,a}$ to be the fiber product of $E_{X,a}\to Y \leftarrow \tilde Y$. Let $\theta_a$ be the pullback of $\theta_{\delta}$ under $\tilde E_{X,a}\to \tilde Y$. 
 
 Note that $q^{\Z}\times q^{\Z}$ acts on $\tilde E_{X,a}$. The action of $(q,q)\in q^{\Z}\times q^{\Z}$ does not change $\theta_a$ and the pullback by the action of $(1,q)$ changes $\theta_a$ to $-{\tilde a}t^{-1}\theta_a$, where $\tilde a$ (resp.\ $t$) is the composition of $\tilde E_{X,a}\to \tilde Y$ and the first (resp.\ second) projection $ \tilde Y\to \mathbb{G}_{m,\log}$.

 \end{para}

\begin{para} Since the class of  $\theta_a$  in $\cQ^{\gp}/M^{\gp}$ is invariant under the action of $q^{\Z}\times q^{\Z}$,  the $M^{\gp}$-torsor associated to this class descends to an $M^{\gp}$-torsor on $E_{X,a}$. We denote it by $\cT(a)$. 
\end{para}

\begin{prop}\label{theta6} Let $a_i\in E^{(q)}(X)$, $n(i)\in \Z$ for $1\leq i\leq r$, and assume that $\sum_{i=1}^r n(i)=0$ and $\prod_{i=1}^r a_i^{n(i)}=1$ in the group $E^{(q)}(X)$. Let $E_{X,a_1, \dots,a_r}$ be the fiber product of $E_{X,a_i}$ ($1\leq i\leq r$) over $X\times E$ which is a log modification of $X \times E$. 

Then on $E_{X,a_1, \dots, a_r}$, $\prod_{i=1}^r \cT(a_i)^{n(i)}$ is a trivial torsor.

\end{prop}

\begin{pf} For $r\geq 1$, let $Y^{(r)}$ (resp.\ $\tilde Y^{(r)}$) be the log modification of $E^{r+1}$ (resp.\ $\tilde E^{r+1}$) defined by the condition on $(t_1, \dots, t_{r+1})$ that $(t_i, t_{r+1})$ belongs to $Y$ (resp.\ $\tilde Y$) for $1\leq i\leq r$. Hence $Y^{(1)}=Y$ (resp.\ $\tilde Y^{(1)}=\tilde Y$). Then $E_{X,a_1, \dots, a_r}$ 
 is the fiber product of $X \times E \to (E^{(q)})^r \times E \leftarrow Y^{(r)}$,
where the first arrow is $(x, t)\mapsto (a_1(x), \dots, a_r(x), t)$. 
Let $\tilde E_{X, a_1, \dots, a_r}$ be the fiber product of $E_{X, a_1, \dots, a_r}\to Y^{(r)}\leftarrow \tilde Y^{(r)}$. We pullback $\theta_{a_i}$ ($1\leq i\leq r$) to $\tilde E_{X,a_1, \dots, a_r}$. Consider $f=\prod_{i=1}^r \theta_{a_i}^{n(i)}$ there. 
For the action of $(q^{\Z})^{r+1}$, $f$ does not change under the action of the element $\gamma_i$ ($1\leq i\leq r$) whose $i$-th component and $r+1$-component are $q$ and whose other components are $1$. For the pullback by the action of $\gamma$ whose $r+1$-th component is $q$ and whose other components are $1$, $\gamma^*(f)f^{-1}\in \mathbb{G}_{m,\log}^{(q)}$ and its image in $E^{(q)}$ is trivial by the assumption on $a_i$. Hence $\gamma^*(f)f^{-1}$ is a locally constant function on $\tilde E_{X,a_1, \dots, a_r}$ with values in $q^{\Z}$, which we denote by $c$. Let $g$ be the section $t^{-c} \prod_{i=1}^r \theta_{a(i)}^{n(i)}$ of $\cQ^{\gp}$ of $\tilde E_{X,a_1, \dots, a_r}$. Then $g$ is invariant under the action of $(q^{\Z})^{r+1}$ and hence is a section of $\cQ^{\gp}$ of $E_{X, a_1, \dots, a_r}$. This gives a section of the torsor $\prod_{i=1}^r \cT(a_i)^{n(i)}$. 
\end{pf}

\begin{lem}\label{KlogM} Let $Y$ be an object of $\frak S$ and let $a,b\in \text{{\rm Pic}}(Y)$. If the images of $a$ and $b$ in $H^1(Y, M^{\gp}_Y)$ coincide, their images in $\gr^1K_0^{\log}(Y)$ coincide. 

\end{lem}

\begin{pf}
This follows from the exact sequence $H^0(Y, M_Y^{\gp}/\cO_Y^\times) \to H^1(Y, \cO_Y^\times)\to H^1(Y, M^{\gp}_Y)$. 
\end{pf}

\begin{para}\label{theta4} Now  (2) in \ref{pf4} follows from
 Lemma \ref{theta8}, Proposition \ref{theta6}, and Lemma \ref{KlogM}.

\end{para}

\noindent {\rm Kazuya KATO
\\
Department of Mathematics
\\
University of Chicago
\\
Chicago, Illinois, 60637, USA}
\\
{\tt kkato@uchicago.edu}

\bigskip

\noindent {\rm Chikara NAKAYAMA
\\
Department of Economics 
\\
Hitotsubashi University 
\\
2-1 Naka, Kunitachi, Tokyo 186-8601, Japan}
\\
{\tt c.nakayama@r.hit-u.ac.jp}

\bigskip

\noindent
{\rm Sampei USUI${}^{\sharp}$
\\
Graduate School of Science
\\
Osaka University}
\\
\footnote[0]{${\sharp}$ \ \ the deceased}
\end{document}